\documentclass{article}
\usepackage{amssymb,latexsym,amsmath,amsthm,mathrsfs}
\usepackage{graphicx}
\newcommand{\comment}[1]{}

\begin{document}
\title{A double demonstration of a theorem of Newton, which gives
a relation between the coefficient of an algebraic equation and the sums
of the powers of its roots\footnote{Presented to the
Berlin Academy on January 12, 1747.
Originally published as
{\em Demonstratio gemina theorematis Neutoniani, quo traditur relatio inter coefficientes
cuiusvis aequationis algebraicae et summas potestatum radicum eiusdem},
Opuscula varii argumenti \textbf{2} (1750), 108--120.
E153 in the Enestr{\"o}m index.
Translated from the Latin by Jordan Bell, School of Mathematics and
Statistics, Carleton University, Ottawa, Canada.
Email: jbell3@connect.carleton.ca}}
\author{Leonhard Euler}
\date{}
\maketitle

1. After an algebraic equation has been given, either from fractions
or by irrationals, and reduced to this form
\[
x^n-Ax^{n-1}+Bx^{n-2}-Cx^{n-3}+Dx^{n-4}-Ex^{n-5}+\cdots \pm N=0,
\]
one demonstrates in analysis that the equation always has as many roots,
either real or imaginary, as unities are contained in 
the exponent $n$ of the highest power. Then indeed it is no less certain
that if all the roots of this equation are $\alpha,\beta,\gamma,\delta,\epsilon,\ldots,\nu$, the coefficients of the terms of the equation $A,B,C,D,E$ etc.
are comprised thus from them
\[
\begin{array}{lllll}
A&=&\textrm{sum of all the roots}&=&\alpha+\beta+\gamma+\delta+\cdots+\nu,\\
B&=&\textrm{sum of the products of two}&=&\alpha\beta+\alpha\gamma+\alpha\delta+\beta\gamma+\textrm{etc.},\\
C&=&\textrm{sum of the products of three}&=&\alpha\beta\gamma+\textrm{etc.},\\
D&=&\textrm{sum of the products of four}&=&\alpha\beta\gamma\delta+\textrm{etc.},\\
E&=&\textrm{sum of the products of five}&=&\alpha\beta\gamma\delta\epsilon+\textrm{etc.}\\
&&\textrm{etc.},&&
\end{array}
\]
and finally the last term $N$ is the product of all the roots $\alpha\beta\gamma\delta\cdots\nu$.

2. Now it will help me to simply and briefly explain this theorem,
whose demonstration I give in this paper,
if I first introduce some notation. Let $\int \alpha$ designate the sum of all
the roots,
$\int \alpha^2$ the sum of the squares of the roots,
$\int \alpha^3$ the sum of the cubes of the roots,
$\int \alpha^4$ the sum of the biquadrates of the roots, and so on,
so that it will be
\[
\begin{array}{lllllllll}
\int \alpha&=&\alpha&+\beta&+\gamma&+\delta&+\epsilon&+\cdots&+\nu,\\
\int \alpha^2&=&\alpha^2&+\beta^2&+\gamma^2&+\delta^2&+\epsilon^2&+\cdots&+\nu^2,\\
\int \alpha^3&=&\alpha^3&+\beta^3&+\gamma^3&+\delta^3&+\epsilon^3&+\cdots&+\nu^3,\\
\int \alpha^4&=&\alpha^4&+\beta^4&+\gamma^4&+\delta^4&+\epsilon^4&+\cdots&+\nu^4,\\
\int \alpha^5&=&\alpha^5&+\beta^5&+\gamma^5&+\delta^5&+\epsilon^5&+\cdots&+\nu^5\\
&&\textrm{etc.}&&&&&&
\end{array}
\]

3. With the rule for this sign explained, Newton affirms that these
sums of powers, which are formed from all the roots, can thus 
be defined by the coefficients of the equation $A,B,C,D,E$ etc.

\begin{eqnarray*}
\int \alpha&=&A,\\
\int \alpha^2&=&A\int \alpha-2B,\\
\int \alpha^3&=&A\int \alpha^2-B\int \alpha+3C,\\
\int \alpha^4&=&A\int \alpha^3-B\int \alpha^2+C\int \alpha-4D,\\
\int \alpha^5&=&A\int \alpha^4-B\int \alpha^3+C\int \alpha^2-D\int \alpha+5E,\\
\int \alpha^6&=&A\int \alpha^5-B\int \alpha^4+C\int \alpha^3-D\int \alpha^2+E\int \alpha-6F\\
&&\textrm{etc.}
\end{eqnarray*}

Newton does not give any demonstration of this theorem, but rather seems
to conclude its truth by continued examples. First, as no demonstration is required that $\alpha =A$, and since
\[
A^2=\alpha^2+\beta^2+\gamma^2+\delta^2+\epsilon^2+\textrm{etc.}+2\alpha\beta+2\alpha\gamma+2\alpha\delta+2\beta\gamma+2\beta\delta+\textrm{etc.},
\]
it will be
\[
A^2=\int \alpha^2+2B
\]
and hence
\[
\alpha^2=A^2-2B=A\alpha-2B;
\]
in a similar way the truth of the following formulae can be shown, but continuing
will take more work.

4. Though many have already shown the truth of this most useful theorem,
 their demonstrations are often based on rules of combinations, which, though
true, depend however on induction. I will thus offer here a double
demonstration,
the first of which completely avoids induction. Indeed the one from the analysis
of the infinite is small, and though it may seem quite remote, everything
comes together perfectly;
however since it can reasonably be objected that the truth of this theorem
should be obtained before the analysis of the infinite is come to,
I add another demonstration, in which nothing is assumed beyond that which
follows immediately from the nature of equations.

\begin{center}
{\Large Demonstration 1.}
\end{center}

5. Let us put
\[
x^n-Ax^{n-1}+Bx^{n-2}-Cx^{n-3}+\cdots \pm N=Z,
\]
and if the roots of the equation $Z=0$, or values of $x$, are
\[
\alpha,\beta,\gamma,\delta,\ldots,\nu,
\]
whose number is $=n$, it will be from the nature of equations
\[
Z=(x-\alpha)(x-\beta)(x-\gamma)(x-\delta)\cdots(x-\nu)
\]
and by taking the logarithm we will obtain
\[
lZ=l(x-\alpha)+l(x-\beta)+l(x-\gamma)+l(x-\delta)+\cdots+l(x-\nu).
\]

But if the differential of these formulae are now taken, it will be
\[
\frac{dZ}{Z}=\frac{dx}{x-\alpha}+\frac{dx}{x-\beta}+\frac{dx}{x-\gamma}+\frac{dx}{x-\delta}+\cdots+\frac{dx}{x-\nu}
\]
and so by dividing by $dx$ it will become
\[
\frac{dZ}{Zdx}=\frac{1}{x-\alpha}+\frac{1}{x-\beta}+\frac{1}{x-\gamma}+\frac{1}{x-\delta}+\cdots+\frac{1}{x-\nu}.
\]
Now let each of these fractions be converted into infinite geometric series in the usual way, as
\begin{eqnarray*}
\frac{1}{x-\alpha}&=&\frac{1}{x}+\frac{\alpha}{x^2}+\frac{\alpha^2}{x^3}+\frac{\alpha^3}{x^4}+\frac{\alpha^4}{x^5}+\frac{\alpha^5}{x^6}+\textrm{etc.},\\
\frac{1}{x-\beta}&=&\frac{1}{x}+\frac{\beta}{x^2}+\frac{\beta^2}{x^3}+\frac{\beta^3}{x^4}+\frac{\beta^4}{x^5}+\frac{\beta^5}{x^6}+\textrm{etc.},\\
\frac{1}{x-\gamma}&=&\frac{1}{x}+\frac{\gamma}{x^2}+\frac{\gamma^2}{x^3}+\frac{\gamma^3}{x^4}+\frac{\gamma^4}{x^5}+\frac{\gamma^5}{x^6}+\textrm{etc.}\\
&&\textrm{etc.}\\
\frac{1}{x-\nu}&=&\frac{1}{x}+\frac{\nu}{x^2}+\frac{\nu^2}{x^3}
+\frac{\nu^3}{x^4}+\frac{\nu^4}{x^5}+\frac{\nu^5}{x^6}+\textrm{etc.}
\end{eqnarray*}

Therefore by collecting these series together, using the signs explained above
$\int \alpha,\int \alpha^2,\int \alpha^3$ etc., because the number of these
series is $=n$,
\[
\frac{dZ}{Zdx}=\frac{n}{x}+\frac{1}{x^2}\int \alpha+\frac{1}{x^3}\int \alpha^2+\frac{1}{x^4}\int \alpha^3+\frac{1}{x^5}\int \alpha^4+\textrm{etc.}
\]

6. On the other hand, since we have put
\[
Z=x^n-Ax^{n-1}+Bx^{n-2}-Cx^{n-3}+Dx^{n-4}-\cdots \pm N,
\]
it will similarly be, by taking differentials,
\[
\frac{dZ}{dx}=nx^{n-1}-(n-1)Ax^{n-2}+(n-2)Bx^{n-3}-(n-3)Cx^{n-4}+(n-4)Dx^{n-5}-\textrm{etc.}
\]
and then let the above formula be put in the same form as $\frac{dZ}{Zdx}$, so that it would be
{\small
\[
\frac{dZ}{Zdx}=\frac{nx^{n-1}-(n-1)Ax^{n-2}+(n-2)Bx^{n-3}-(n-3)Cx^{n-4}+(n-4)Dx^{n-5}-\textrm{etc.}}{x^n-Ax^{n-1}+Bx^{n-2}-Cx^{n-3}+Dx^{n-4}-\textrm{etc.}},
\]
}
and this fraction will therefore be equal to the series found previously
\[
\frac{n}{x}+\frac{1}{x^2}\int \alpha+\frac{1}{x^3}\int \alpha^2+\frac{1}{x^4}\int \alpha^3+\frac{1}{x^5}\int \alpha^4+\textrm{etc.}
\]
Now if each expression found for $\frac{dZ}{Zdx}$ is multiplied by the one
denominator
\[
x^n-Ax^{n-1}+Bx^{n-2}-Cx^{n-3}+Dx^{n-4}-\textrm{etc.},
\]
this equation will result
\[
\begin{split}
&nx^{n-1}-(n-1)Ax^{n-2}+(n-2)Bx^{n-3}-(n-3)Cx^{n-4}+(n-4)Dx^{n-5}-\textrm{etc.}\\
&=\\
&\begin{array}{llllll}
nx^{n-1}&+x^{n-2}\int \alpha&+x^{n-3}\int \alpha^2&+x^{n-4}\int \alpha^3&+x^{n-5}\int \alpha^4&+\textrm{etc.}\\
&-nAx^{n-2}&-Ax^{n-3}\int \alpha&-Ax^{n-4}\int \alpha^2&-Ax^{n-5}\int \alpha^3&-\textrm{etc.}\\
&&+nBx^{n-3}&+Bx^{n-4}\int \alpha&+Bx^{n-5}\int \alpha^2&+\textrm{etc.}\\
&&&-nCx^{n-4}&-Cx^{n-5}\int \alpha&-\textrm{etc.}\\
&&&&+nDx^{n-5}&+\textrm{etc.}
\end{array}
\end{split}
\]

7. Now just as much as both the first terms $nx^{n-1}$ are equal,
it is also necessary that the second, third, fourth etc. are equal to each
other; thus the following equations arise
\begin{eqnarray*}
-(n-1)A&=&\int \alpha-nA,\\
+(n-2)B&=&\int \alpha^2-A\int \alpha+nB,\\
-(n-3)C&=&\int \alpha^3-A\int \alpha^2+B\int \alpha-nC,\\
+(n-4)D&=&\int \alpha^4-A\int \alpha^3+B\int \alpha^2-C\int \alpha+nD\\
&&\textrm{etc.}
\end{eqnarray*}
and the law by which these equations progress is at immediately clear. 
Further, from these, the formulae can be obtained which comprise the theorem of
Newton, namely
\begin{eqnarray*}
\int \alpha&=&A,\\
\int \alpha^2&=&A\int \alpha-2B,\\
\int \alpha^3&=&A\int \alpha^2-B\int \alpha+3C,\\
\int \alpha^4&=&A\int \alpha^3-B\int \alpha^2+C\int \alpha-4D,\\
\int \alpha^5&=&A\int \alpha^4-B\int \alpha^3+C\int \alpha^2-D\int \alpha+5E\\
&&\textrm{etc.}
\end{eqnarray*}
This is the first demonstration of the proposed theorem.

\begin{center}
{\Large Demonstration 2.}
\end{center} 

8. So this demonstration can be clearly understood, I will fix an equation
of a determinate degree, but one can see that this extends equally
well to any degree. Thus let the given equation be of the fifth degree
\[
x^5-Ax^4+Bx^3-Cx^2-E=0,
\]
all of whose roots are $\alpha,\beta,\gamma,\delta,\epsilon$.
Then because any root substituted in place of $x$ satisfies the equation,
it will be
\[
\begin{array}{lllllll}
\alpha^5&-A\alpha^4&+B\alpha^3&-C\alpha^2&+D\alpha&-E&=0,\\
\beta^5&-A\beta^4&+B\beta^3&-C\beta^2&+D\beta&-E&=0,\\
\gamma^5&-A\gamma^4&+B\gamma^3&-C\gamma^2&+D\gamma&-E&=0,\\
\delta^5&-A\delta^4&+B\delta^3&-C\delta^2&+D\delta&-E&=0,\\
\epsilon^5&-A\epsilon^4&+B\epsilon^3&-C\epsilon^2&+D\epsilon&-E&=0.
\end{array}
\]
These equations can be collected into one sum, and using the above sign (\S 2) one will have
\[
\int \alpha^5-A\int \alpha^4+B\int \alpha^3-C\int \alpha^2+D\int \alpha-5E=0
\]
or
\[
\int \alpha^5=A\int \alpha^4-B\int \alpha^3+C\int \alpha^2-D\int \alpha+5E.
\]

9. It is clear then that if the given equation were of any degree
\[
x^n-Ax^{n-1}+Bx^{n-2}-Cx^{n-3}+Dx^{n-4}-\cdots \pm Mx \mp N=0,
\] 
where in the last term, the top of the ambiguous signs occurs if the exponent
$n$ is an odd number, and the bottom if even, it will likewise be
\[
\int \alpha^n=A\int \alpha^{n-1}-B\int \alpha^{n-2} +C \int \alpha^{n-3}-\cdots \mp M\int \alpha \pm nN,
\]
where $\alpha$ indicates any root of this equation. Thus the truth of
the Theorem of Newton is now shown for this one case. It therefore
remains for us to demonstrate the truth of it for both the higher and for the lower
 powers of the roots.

10. Indeed for the higher powers the same argument applies; for if
the values $\alpha,\beta,\gamma,\delta,\epsilon$ satisfy the equation
\[
x^5-Ax^4+Bx^3-Cx^2+Dx-E=0,
\]
then they also will satisfy the following equations
\begin{eqnarray*}
x^6-Ax^5+Bx^4-Cx^3+Dx^2-Ex&=&0,\\
x^7-Ax^6+Bx^5-Cx^4+Dx^3-Ex^2&=&0,\\
x^8-Ax^7+Bx^6-Cx^5+Dx^4-Ex^3&=&0,\\
\textrm{etc.}&&
\end{eqnarray*}
Therefore also if in each equation all the values $\alpha,\beta,\gamma,\delta,\epsilon$ are substituted for $x$ and the aggregate collected, it will be
\[
\begin{split}
&\int \alpha^6=A\int \alpha^5-B\int \alpha^4+C\int \alpha^3-D\int \alpha^2+E\int \alpha,\\
&\int \alpha^7=A\int \alpha^6-B\int \alpha^5+C\int \alpha^4-D\int \alpha^3+E\int \alpha^2,\\
&\int \alpha^8=A\int \alpha^7+C\int \alpha^5-D\int \alpha^4+E\int \alpha^3\\
&\textrm{etc.}
\end{split}
\]

11. Thus if $\alpha$ denotes any root of this equation
\[
x^n-Ax^{n-1}+Bx^{n-2}-Cx^{n-3}+Dx^{n-4}-\cdots \pm Mx \mp N=0,
\]
it will not only be, as we have already found,
\[
\int \alpha^n=A\int \alpha^{n-1}-B\int \alpha^{n-2}+C\int \alpha^{n-3}-D\int \alpha^{n-4} + \cdots \mp M \int \alpha \pm nN,
\]
but by proceeding to higher powers it will also be
{\small 
\[
\begin{split}
&\int \alpha^{n+1}=A\int \alpha^n-B\int \alpha^{n-1}+C\int \alpha^{n-2}-D\int \alpha^{n-3}+\cdots \mp M\int \alpha^2 \pm N\int \alpha,\\
&\int \alpha^{n+2}=A\int \alpha^{n+1}-B\int \alpha^n+C\int \alpha^{n-1}-D\int \alpha^{n-2}+\cdots \mp M\int \alpha^3 \pm N \int \alpha^2,\\
&\int \alpha^{n+3}=A\int \alpha^{n+2}-B\int \alpha^{n+1}+C\int \alpha^n-D\int \alpha^{n-1}+\cdots \mp M \int \alpha^4 \pm N\int \alpha^3,\\
&\textrm{etc.}
\end{split}
\]
}
and indeed in general, if any number $m$ is added to $n$, it will be
{\small
\[
\int \alpha^{n+m}=A\int \alpha^{n+m-1}-B\int \alpha^{n+m-2}+C\int \alpha^{n+m-3} -\cdots \mp M \int \alpha^{m+1} \pm N \int \alpha^m.
\]
}
It should be noted here that if $m=0$, because all the powers
$\alpha^0=1,\beta^0=1,\gamma^0=1$ etc. and the number of these letters $=n$, then $\int \alpha^0=n$,
and in this case the first formula that was found is included in this expression.

12. This expression could also be used if a negative number is taken
for $m$, and then for the assumed equation of the fifth degree
\[
x^5-Ax^4+Bx^3-Cx^2+Dx-E=0
\]
the following formulae would respectively occur
\begin{eqnarray*}
\int \alpha^4&=&A\int \alpha^3-B\int \alpha^2+C\int \alpha^1-D\int \alpha^0+E\int \alpha^{-1},\\
\int \alpha^3&=&A\int \alpha^2-B\int \alpha^1+C\int \alpha^0-D\int \alpha^{-1}+E\int \alpha^{-2},\\
\int \alpha^2&=&A\int \alpha^1-B\int \alpha^0+C\int \alpha^{-1}-D\int \alpha^{-2}+E\int \alpha^{-3}\\
&&\textrm{etc.},
\end{eqnarray*}
but these formulae are different from those which the theorem contains.
For it needs to be demonstrated that
\begin{eqnarray*}
\int \alpha^4&=&A\int \alpha^3-B\int \alpha^2+C\int \alpha-4D,\\
\int \alpha^3&=&A\int \alpha^2-B\int \alpha+3C,\\
\int \alpha^2&=&A\int \alpha-2B,\\
\int \alpha&=&A.
\end{eqnarray*}
I will therefore show the truth of these formulae
in the following way.

13. Namely, given an equation of the fifth degree
\[
x^5-Ax^4+Bx^3-Cx^2+Dx-E=0
\]
let us form the following equations of lower degrees by using its coefficients

\begin{tabular}{rll}
I.&$x-A=0$.&Let $p$ be its root.\\
II.&$x^2-Ax+B=0$.&Let $q$ be any root.\\
III.&$x^3-Ax^2+Bx-C=0$.&Let $r$ be any root.\\
IV.&$x^4-Ax^3+Bx^2-Cx+D=0$.&Let $s$ be any root.
\end{tabular}

The roots of these equations, even though they are different from
each other, constitute however for all these equations the same sum
$=A$. Next, with the first removed, the sum of the products from two
roots will everywhere be the same $=B$. Then
the sum of the products from three roots will everywhere $=C$,
besides of course equations I and II, where $C$ does not occur. 
Similarly in the given equation IV etc. the sum of the products
from four roots will be the same $=D$.

14. Therefore in the equations of lower degrees, whose roots are denoted
respectively by the letters $p,q,r,s$, if the letter $\alpha$ denotes an arbitrary
root of the given equation,  it will be
\[
\begin{split}
&\int \alpha=\int s=\int r=\int q=\int p,\\
&\int \alpha^2=\int s^2=\int r^2=\int q^2,\\
&\int \alpha^3=\int s^3=\int r^3,\\
&\int \alpha^4=\int s^4.
\end{split}
\]
But on the other hand, it follows, from what we demonstrated above in \S 9, that
\begin{eqnarray*}
\int p&=&A,\\
\int q^2&=&A\int q-2B,\\
\int r^3&=&A\int r^2-B\int r+3C,\\
\int s^4&=&A\int s^3-B\int s^2+C\int s-4D.
\end{eqnarray*} 
Now therefore for the given equation of the fifth degree
\[
x^5-Ax^4+Bx^3-Cx^2+Dx-E=0
\]
these formulae arise
\begin{eqnarray*}
\int \alpha&=&A,\\
\int \alpha^2&=&A\int \alpha-2B,\\
\int \alpha^3&=&A\int \alpha^2-B\int \alpha+3C,\\
\int \alpha^4&=&A\int \alpha^3-B\int \alpha^2+C\int \alpha-4D.
\end{eqnarray*}

16. Therefore in a given equation of any degree
\[
x^n-Ax^{n-1}+Bx^{n-2}-Cx^{n-3}+Dx^{n-4}-\cdots \pm N=0
\]
if the letter $\alpha$ indicates an arbitrary root, it will be
\begin{eqnarray*}
\int \alpha&=&A,\\
\int \alpha^2&=&A\int \alpha-2B,\\
\int \alpha^3&=&A\int \alpha^2-B\int \alpha+3C,\\
\int \alpha^4&=&A\int \alpha^3-B\int \alpha^2+C\int \alpha-4D,\\
\int \alpha^5&=&A\int \alpha^4-B\int \alpha^3+C\int \alpha^2-D\int \alpha+5E\\
&&\textrm{etc.}
\end{eqnarray*}
In this way the truth of the Theorem of Newton is also demonstrated. 

\end{document}